\documentclass{amsart}

\usepackage{amsfonts,amssymb,verbatim,amsmath,amsthm,latexsym,textcomp,amscd}
\usepackage{latexsym,amsfonts,amssymb,epsfig,verbatim}
\usepackage{amsmath,amsthm,amssymb,latexsym,graphics,textcomp}
\usepackage{graphicx}
\usepackage{color}
\usepackage{url}

\input xy
\xyoption{all}

\begin{document}

\newtheorem{theorem}{Theorem}[section]
\newtheorem{prop}[theorem]{Proposition}
\newtheorem{lemma}[theorem]{Lemma}
\newtheorem{cor}[theorem]{Corollary}
\newtheorem{definition}[theorem]{Definition}
\newtheorem{conj}[theorem]{Conjecture}
\newtheorem{rmk}[theorem]{Remark}
\newtheorem{claim}[theorem]{Claim}
\newtheorem{defth}[theorem]{Definition-Theorem}

\newcommand{\boundary}{\partial}
\newcommand{\C}{{\mathbb C}}
\newcommand{\integers}{{\mathbb Z}}
\newcommand{\natls}{{\mathbb N}}
\newcommand{\ratls}{{\mathbb Q}}
\newcommand{\bbR}{{\mathbb R}}
\newcommand{\proj}{{\mathbb P}}
\newcommand{\lhp}{{\mathbb L}}
\newcommand{\tube}{{\mathbb T}}
\newcommand{\cusp}{{\mathbb P}}
\newcommand\AAA{{\mathcal A}}
\newcommand\BB{{\mathcal B}}
\newcommand\CC{{\mathcal C}}
\newcommand\DD{{\mathcal D}}
\newcommand\EE{{\mathcal E}}
\newcommand\FF{{\mathcal F}}
\newcommand\GG{{\mathcal G}}
\newcommand\HH{{\mathcal H}}
\newcommand\II{{\mathcal I}}
\newcommand\JJ{{\mathcal J}}
\newcommand\KK{{\mathcal K}}
\newcommand\LL{{\mathcal L}}
\newcommand\MM{{\mathcal M}}
\newcommand\NN{{\mathcal N}}
\newcommand\OO{{\mathcal O}}
\newcommand\PP{{\mathcal P}}
\newcommand\QQ{{\mathcal Q}}
\newcommand\RR{{\mathcal R}}
\newcommand\SSS{{\mathcal S}}
\newcommand\TT{{\mathcal T}}
\newcommand\UU{{\mathcal U}}
\newcommand\VV{{\mathcal V}}
\newcommand\WW{{\mathcal W}}
\newcommand\XX{{\mathcal X}}
\newcommand\YY{{\mathcal Y}}
\newcommand\ZZ{{\mathcal Z}}
\newcommand\CH{{\CC\HH}}
\newcommand\PEY{{\PP\EE\YY}}
\newcommand\MF{{\MM\FF}}
\newcommand\RCT{{{\mathcal R}_{CT}}}
\newcommand\PMF{{\PP\kern-2pt\MM\FF}}
\newcommand\FL{{\FF\LL}}
\newcommand\PML{{\PP\kern-2pt\MM\LL}}
\newcommand\GL{{\GG\LL}}
\newcommand\Pol{{\mathcal P}}
\newcommand\half{{\textstyle{\frac12}}}
\newcommand\Half{{\frac12}}
\newcommand\Mod{\operatorname{Mod}}
\newcommand\Area{\operatorname{Area}}
\newcommand\ep{\epsilon}
\newcommand\hhat{\widehat}
\newcommand\Proj{{\mathbf P}}
\newcommand\U{{\mathbf U}}
 \newcommand\Hyp{{\mathbf H}}
\newcommand\D{{\mathbf D}}
\newcommand\Z{{\mathbb Z}}
\newcommand\R{{\mathbb R}}
\newcommand\Q{{\mathbb Q}}
\newcommand\E{{\mathbb E}}
\newcommand\til{\widetilde}
\newcommand\length{\operatorname{length}}
\newcommand\tr{\operatorname{tr}}
\newcommand\gesim{\succ}
\newcommand\lesim{\prec}
\newcommand\simle{\lesim}
\newcommand\simge{\gesim}
\newcommand{\simmult}{\asymp}
\newcommand{\simadd}{\mathrel{\overset{\text{\tiny $+$}}{\sim}}}
\newcommand{\ssm}{\setminus}
\newcommand{\diam}{\operatorname{diam}}
\newcommand{\pair}[1]{\langle #1\rangle}
\newcommand{\T}{{\mathbf T}}
\newcommand{\inj}{\operatorname{inj}}
\newcommand{\pleat}{\operatorname{\mathbf{pleat}}}
\newcommand{\short}{\operatorname{\mathbf{short}}}
\newcommand{\vertices}{\operatorname{vert}}
\newcommand{\collar}{\operatorname{\mathbf{collar}}}
\newcommand{\bcollar}{\operatorname{\overline{\mathbf{collar}}}}
\newcommand{\I}{{\mathbf I}}
\newcommand{\tprec}{\prec_t}
\newcommand{\fprec}{\prec_f}
\newcommand{\bprec}{\prec_b}
\newcommand{\pprec}{\prec_p}
\newcommand{\ppreceq}{\preceq_p}
\newcommand{\sprec}{\prec_s}
\newcommand{\cpreceq}{\preceq_c}
\newcommand{\cprec}{\prec_c}
\newcommand{\topprec}{\prec_{\rm top}}
\newcommand{\Topprec}{\prec_{\rm TOP}}
\newcommand{\fsub}{\mathrel{\scriptstyle\searrow}}
\newcommand{\bsub}{\mathrel{\scriptstyle\swarrow}}
\newcommand{\fsubd}{\mathrel{{\scriptstyle\searrow}\kern-1ex^d\kern0.5ex}}
\newcommand{\bsubd}{\mathrel{{\scriptstyle\swarrow}\kern-1.6ex^d\kern0.8ex}}
\newcommand{\fsubeq}{\mathrel{\raise-.7ex\hbox{$\overset{\searrow}{=}$}}}
\newcommand{\bsubeq}{\mathrel{\raise-.7ex\hbox{$\overset{\swarrow}{=}$}}}
\newcommand{\tw}{\operatorname{tw}}
\newcommand{\base}{\operatorname{base}}
\newcommand{\trans}{\operatorname{trans}}
\newcommand{\rest}{|_}
\newcommand{\bbar}{\overline}
\newcommand{\UML}{\operatorname{\UU\MM\LL}}
\newcommand{\EL}{\mathcal{EL}}
\newcommand{\tsum}{\sideset{}{'}\sum}
\newcommand{\tsh}[1]{\left\{\kern-.9ex\left\{#1\right\}\kern-.9ex\right\}}
\newcommand{\Tsh}[2]{\tsh{#2}_{#1}}
\newcommand{\qeq}{\mathrel{\approx}}
\newcommand{\Qeq}[1]{\mathrel{\approx_{#1}}}
\newcommand{\qle}{\lesssim}
\newcommand{\Qle}[1]{\mathrel{\lesssim_{#1}}}
\newcommand{\simp}{\operatorname{simp}}
\newcommand{\vsucc}{\operatorname{succ}}
\newcommand{\vpred}{\operatorname{pred}}
\newcommand\fhalf[1]{\overrightarrow {#1}}
\newcommand\bhalf[1]{\overleftarrow {#1}}
\newcommand\sleft{_{\text{left}}}
\newcommand\sright{_{\text{right}}}
\newcommand\sbtop{_{\text{top}}}
\newcommand\sbot{_{\text{bot}}}
\newcommand\sll{_{\mathbf l}}
\newcommand\srr{_{\mathbf r}}
\newcommand\geod{\operatorname{\mathbf g}}
\newcommand\vol{\operatorname{\mathrm vol}}
\newcommand\mtorus[1]{\boundary U(#1)}
\newcommand\A{\mathbf A}
\newcommand\Aleft[1]{\A\sleft(#1)}
\newcommand\Aright[1]{\A\sright(#1)}
\newcommand\Atop[1]{\A\sbtop(#1)}
\newcommand\Abot[1]{\A\sbot(#1)}
\newcommand\boundvert{{\boundary_{||}}}
\newcommand\storus[1]{U(#1)}
\newcommand\Momega{\omega_M}
\newcommand\nomega{\omega_\nu}
\newcommand\twist{\operatorname{tw}}
\newcommand\modl{M_\nu}
\newcommand\MT{{\mathbb T}}
\newcommand\Teich{{\mathcal T}}
\renewcommand{\Re}{\operatorname{Re}}
\renewcommand{\Im}{\operatorname{Im}}

\title{On Sullivan's construction of eigenfunctions via exit times of Brownian motion}

\author{Kingshook Biswas}
\address{Indian Statistical Institute, Kolkata, India. Email: kingshook@isical.ac.in}

\begin{abstract} The purpose of this note is to give details for an argument of Sullivan to construct eigenfunctions of the Laplacian
on a Riemannian manifold using exit times of Brownian motion \cite{sullivanpos}.
Let $X$ be a complete, simply connected Riemannian manifold of pinched negative sectional curvature.
Let $\lambda_1 = \lambda_1(X) < 0$ be the supremum of the spectrum of the Laplacian on $L^2(X)$, and let
$D \subset X$ be a bounded domain in $X$ with smooth boundary. Let $(B_t)_{t \geq 0}$ be Brownian motion on $X$ and let
$\tau = \tau_D$ be
the first exit time of Brownian motion from $D$. For each
$\lambda \in \mathbb{C}$ with $\hbox{Re } \ \lambda > \lambda_1$ and $x \in D$, we show that
for any continuous function $\phi : \partial D \to \mathbb{C}$, the function
$$
h(x) = \mathbb{E}_x(e^{-\lambda \tau} \phi(B_{\tau})) \ , \ x \in D,
$$
is an eigenfunction of the Laplacian on $D$ with eigenvalue $\lambda$ and boundary value $\phi$.
\end{abstract}

\bigskip

\maketitle

\tableofcontents

\section{Introduction}

\medskip

The purpose of this note is to give the details for a construction of Sullivan of eigenfunctions of the Laplacian on compact domains with smooth boundary using
exit times for the Brownian motion which is just briefly sketched in \cite{sullivanpos}. For simplicity,
we restrict ourselves to manifolds with pinched negative curvature, though the arguments should still hold in a more
general setting where the manifold is stochastically complete and the Brownian motion is transient.

\medskip

Sullivan's result, which is just stated in \cite{sullivanpos} without any proof or reference unfortunately, is the following:

\medskip

\begin{theorem} \label{mainthm} Let $X$ be a complete, simply connected Riemannian manifold of pinched negative sectional curvature,
and let $(B_t)_{t \geq 0}$ be Brownian motion on $X$. Let $\lambda_1(X) < 0$ be the supremum of the Laplacian $\Delta$ on $L^2(X)$.
Let $D \subset X$ be a precompact domain in $X$ with smooth boundary, and let $\tau$ be the first exit time from $D$ for Brownian motion.
Then for any $\lambda \in \mathbb{C}$ with $\hbox{Re } \ \lambda > \lambda_1$, and for any continuous function $\phi : \partial D \to \mathbb{C}$, the function
$$
h(x) := \mathbb{E}_x(e^{-\lambda \tau} \phi(B_{\tau}))
$$
is $C^{\infty}$ on $D$ and is an eigenfunction of $\Delta$ on $D$ with eigenvalue $\lambda$ and boundary value $\phi$, meaning $\Delta h = \lambda h$ and
$$
h(x) \to \phi(\xi)
$$
as $x \in D \to \xi \in \partial D$.
\end{theorem}

\medskip

For $\lambda = 0$, this is just the well-known classical solution of the Dirichlet problem using Brownian motion. The
above theorem for other values of $\lambda > \lambda_1$ must be well known to experts, but owing to the lack of a proof
or reference in Sullivan, it seemed worthwhile
to write down the details. The article is organized as follows: in section 2 we present basic facts about the heat kernel
on a complete Riemannian manifold. In section 3 we describe the construction of Brownian motion from the heat semigroup.
In section 4 we show that the infinitesimal generator of the heat semigroup restricted to $C^{\infty}_c(X)$ is given by the
Laplacian. In section 5 we introduce the Dirichlet heat kernel of a bounded domain, while in section 6 we give the proof of the above theorem.

\medskip

\section{The heat kernel}

\medskip

Let $X$ be a complete, simply connected Riemannian manifold with pinched negative sectional curvature $-b^2 \leq K \leq -a^2$.
Let $\Delta$ be the Laplacian on $X$, acting on functions $f \in C^{\infty}(X)$ by
$$
\Delta f = \hbox{div}(\nabla f)
$$
Gaffney \cite{gaffney} showed that the densely defined operator $\Delta$ on $L^2(X)$ with domain $C^{\infty}_c(X)$
is essentially self-adjoint, hence it has
a unique self-adjoint extension on $L^2(X)$, also denoted by $\Delta$. The domain of $\Delta$ is given by $f$ in $L^2(X)$ such that
$\Delta f$ (in the sense of distributions) is in $L^2(X)$. Define $\lambda_1 = \lambda_1(X) \leq 0$ by
$$
\lambda_1(X) := - \inf_{\phi} \frac{\int_X |\nabla \phi|^2 dvol}{\int_X |\phi|^2 dvol} \ ,
$$
where the infimum is taken over all non-zero $\phi$ in $C^{\infty}_c(X)$.
Then the spectrum of $\Delta$ on $L^2(X)$ is contained in $(-\infty, \lambda_1]$ and $\lambda_1$ is the supremum of the spectrum \cite{chavel}.
The functional calculus for self-adjoint operators then allows us to define the {\it heat semigroup} $(e^{t\Delta})_{t \geq 0}$
as a semigroup of bounded operators on $L^2(X)$ satisfying
$$
e^{(s+t)\Delta} = e^{s\Delta} e^{t\Delta} \ , s,t \geq 0.
$$
Mckean \cite{mckean}
showed that the upper bound on sectional curvature $K \leq -a^2$ implies
that
$$
\lambda_1 \leq -(n-1)^2 a^2/4 \ ,
$$
where $n$ is the dimension of $X$, in particular $\lambda_1 < 0$.

\medskip

It is known that the action of the semigroup $e^{t\Delta}$ on $L^2(X)$ is given by integrating against
a kernel $p(t, x, y)$ called the
{\it heat kernel} \cite{dodziuk}.
The heat kernel is a positive smooth function $p : (0, \infty) \times X \times X \to (0, \infty)$ satisfying
\begin{align*}
\frac{\partial p}{\partial t} & = \Delta_y p \\
p(t, x, .) & \to \delta_x \ \hbox{ as } \ t \to 0 \\
\end{align*}
where the second condition means
$$
\int_X p(t, x, y) f(y) dvol(y) \to f(x) \ \hbox{ as } t \to 0
$$
for all bounded continuous functions $f$ (the heat kernel satisfies $\int_X p(t, x, y) dvol(y) \leq 1$ so the above integral is well-defined for $f$ bounded).
The semigroup $e^{t\Delta}$ is given by
$$
e^{t\Delta} f (x) = \int_X p(t,x,y) f(y) dvol(y)
$$
for $f \in L^2(X)$.
For the existence of the heat kernel,
see \cite{chavel}, \cite{dodziuk}. If $X$ has Ricci curvature bounded from below, which is the case with our hypothesis
of sectional curvature bounded below, then the heat kernel is unique \cite{chavel}.

\medskip

Grigoryan \cite{grigoryan} gave upper bounds on the heat kernel for complete, non-compact manifolds in a general setting.
In our case of a complete simply connected manifold of pinched negative curvature, these specialize to the following estimates:

\medskip

There exist constants $c, C > 0$ and $D > 4$ such that for all $x, y \in X, t > 0$, we have
\begin{equation} \label{short time}
p(t,x,y) \leq C t^{-n/2} \exp\left(c\lambda_1 t - \frac{d(x,y)^2}{Dt}\right)
\end{equation}
For long-time asymptotics, this can be improved: for all $x,y \in X, t > 1$, we have
\begin{equation} \label{long time}
p(t,x,y) \leq K \left( 1 + \frac{d(x,y)^2}{t} \right)^{1 + \frac{n}{2}} \exp\left( \lambda_1 t - \frac{d(x,y)^2}{4t} \right)
\end{equation}

\medskip

Since $X$ has sectional curvature bounded below, the volume growth of $X$ is
at most exponential, and there are constants $K, h > 0$ such that the Jacobian $J(x, v, r)$ of the map $v \in T^1_x X \mapsto exp_x(rv)$
at $v$ satisfies a bound
\begin{equation} \label{jacobian}
J(x,v,r) \leq K e^{hr}
\end{equation}
for all $x \in X, v \in T^1_x X, r > 0$. Together with (\ref{short time}) above, this leads to the following lemma:

\medskip

\begin{lemma} \label{outside ball} There are constants $\kappa, \eta > 0$ such that for all $x \in X$, for $R > 0$ and $t > 0$  such
that $R \geq (Dh) t$ we have
$$
\int_{X - B(x, R)} p(t,x,y) dvol(y) \leq \kappa e^{-\eta R^2/t}
$$
\end{lemma}

\medskip

\noindent{\bf Proof:} Integrating in geodesic polar coordinates centered at $x$ and using the estimates (\ref{short time}), (\ref{jacobian}) gives
\begin{align*}
\int_{X - B(x, R)} p(t,x,y) dvol(y) & \leq C t^{-n/2} e^{c\lambda_1 t} \int_{R}^{\infty} e^{-r^2/(Dt)} Ke^{hr} dr \\
                                   & \leq KC t^{-n/2} \sqrt{Dt} \int_{\frac{R}{\sqrt{Dt}} - \frac{\sqrt{Dt}h}{2}}^{\infty} e^{-s^2} ds
\end{align*}
Now $R \geq (Dh)t$ gives $\frac{R}{\sqrt{Dt}} - \frac{\sqrt{Dt}h}{2} \geq \frac{R}{2\sqrt{Dt}}$, then using a
standard bound for the tail of a Gaussian integral
gives, for some constant $C_0 > 0$,
$$
\int_{X - B(x, R)} p(t,x,y) dvol(y) \leq C_0 t^{-(n-1)/2} e^{-R^2/(4Dt)}
$$
and the last expression above is bounded by $\kappa e^{-\eta R^2/t}$, choosing $\kappa$ large enough and $\eta$ small enough. $\diamond$

\medskip

\section{Construction of Brownian motion}

\medskip

We first briefly recall the correspondence between Markov processes and (certain) semigroups of bounded operators, and then
explain how the heat semigroup $(e^{t\Delta})_{t \geq 0}$ can be used to construct Brownian motion on a manifold.

\medskip

\subsection{Markov processes and semigroups}

\medskip

Let $(X, d)$ be a locally compact, separable metric space equipped with its Borel sigma-algebra $\mathcal{B}(X)$.
Let $B(X)$ denote the Banach space of bounded measurable functions on $X$ equipped with the supremum norm and $C_0(X)$
the closed subspace of continuous functions on $X$ vanishing at infinity.

\medskip

We recall that {\it conditional expectation} is defined as follows: given a probability space $(\Omega, \mathcal{F}, \mathbb{P})$,
for any $\phi \in L^1(\Omega, \mathcal{F}, \mathbb{P})$ and any sigma algebra $\mathcal{G} \subset \mathcal{F}$, the conditional expectation
of $\phi$ given $\mathcal{G}$ is the unique $\mathcal{G}$-measurable random variable $\mathbb{E}(\phi|\mathcal{G}) \in L^1(\Omega, \mathcal{G}, \mathbb{P})$
such that
$$
\int_A \mathbb{E}(\phi|\mathcal{G}) d\mathbb{P} = \int_A \phi d\mathbb{P}
$$
for all $A \in \mathcal{G}$ (the existence of the conditional expectation
follows from an application of the Radon-Nikodym theorem to the signed measure
$A \in \mathcal{G} \mapsto \int_A \phi d\mathcal{P}$ on $\mathcal{G}$). For a random variable $Y$ on $\Omega$, we define the conditional expectation
$\mathbb{E}(\phi|Y)$ to be the conditional expectation of $\phi$ given the sigma algebra generated by $Y$.

\medskip

Given a probability measure $\nu$ on $X$, a
{\it Markov process} on $X$ with initial distribution $\nu$ is a collection $(B_t)_{t \geq 0}$ of $X$-valued
random variables on a probability space $(\Omega, \mathcal{F}, \mathbb{P})$ (i.e. each $B_t$ is a measurable map $B_t : \Omega \to X$)
such that for all bounded measurable functions $f$ on $X$ we have the {\it Markov Property}
$$
\mathbb{E}(f(B_t)|\mathcal{F}_s) = \mathbb{E}(f(B_t)|B_s)
$$
for all $t \geq s$  (where $\mathcal{F}_s \subset \mathcal{F}$ is the sub-sigma algebra generated by the maps $B_u, 0 \leq u \leq s$),
and such that $B_0$ has distribution $\nu$, i.e.
$$
\mathbb{P}(B_0 \in E) = \nu(E)
$$
for all Borel sets $E \subset X$.
By a {\it sample path} of the process, we mean a path in $X$ of the form $t \mapsto B_t(\omega)$, for some $\omega \in \Omega$.

\medskip

Let $(T(t))_{t \geq 0}$ be a semigroup of bounded operators on $B(X)$. We say that a Markov process $(B_t)_{t \geq 0}$
on $X$ {\it corresponds} to the semigroup $(T(t))_{t \geq 0}$ if for all $f \in B(X)$ and $s,t \geq 0$ we have
\begin{equation} \label{corresponds}
\mathbb{E}(f(B_{s+t})| \mathcal{F}_s) = (T(t)f)(B_s)
\end{equation}
If the above condition holds, then using the Markov property one can show that the initial distribution $\nu$ together
with the semigroup $(T(t))_{t \geq 0}$ determine the {\it finite dimensional distributions} of the process $(B_t)_{t \geq 0}$, i.e. all the probabilities of the form
$\mathbb{P}(B_{t_1} \in E_1, \dots, B_{t_n} \in E_n)$ are determined, for any finite sequence $0 \leq t_1 < \dots < t_n$ and any Borel sets
$E_1, \dots, E_n \in \mathcal{B}(X)$.

\medskip

Conversely, given a semigroup satisfying certain properties, one can construct a Markov process which corresponds to the semigroup and whose
sample paths have good regularity properties. The precise statement is the following:

\medskip

Let $(T(t))_{t \geq 0}$ be a semigroup of bounded operators on $B(X)$ which satisfies the following properties:

\medskip

\noindent (a) $||T(t)|| \leq 1$ for all $t \geq 0$ ({\it contraction}).

\medskip

\noindent (b) $T(t)f \geq 0$ for all $f \geq 0$ ({\it positivity}).

\medskip

\noindent (c) For all $f \in C_0(X)$, we have $T(t)f \in C_0(X)$ for all $t \geq 0$ ({\it Feller property}).

\medskip

\noindent (d) For all $f \in C_0(X)$, we have $||T(t)f - f|| \to 0$ as $t \to 0$ ({\it strongly continuous}).

\medskip

\noindent (e) $T(t)1 = 1$ for all $t \geq 0$ (where $1$ denotes the function constant equal to 1) ({\it conservative}).

\medskip

A semigroup satisfying these properties is called a {\it Feller semigroup}.

\medskip

A classical result from probability theory (see for e.g. \cite{ethierkurtz}, Theorem 2.7 of Chapter 4) asserts that given a Feller semigroup $(T(t))_{t \geq 0}$
on a locally compact, separable metric space $X$, for any probability measure $\nu$ on $X$ there is a Markov process $(B_t)_{t \geq 0}$ with initial
distribution $\nu$ which corresponds to the semigroup $(T(t))_{t \geq 0}$, such that for $\mathbb{P}$-almost every $\omega$ the sample path $t \mapsto B_t(\omega)$ is
a {\it cadlag} path, i.e. the path is right-continuous with left-limits existing for all $t$.

\medskip

Let $D_X[0, \infty)$ denote the space of cadlag paths on $X$, for $t \geq 0$ let $\pi_t : D_X([0, \infty)) \to X$ be the map defined by
$\pi_t(\gamma) := \gamma(t)$, and let $\mathcal{D}$ be the sigma-algebra on $D_X([0, \infty))$ generated by the maps $\pi_t, t \geq 0$. Given $x \in X$,
by the above theorem, we have a Markov process $(B_t)_{t \geq 0}$ on some $(\Omega, \mathcal{F}, \mathbb{P})$
with initial distribution $\delta_x$ (the Dirac mass at $x$)
which corresponds to the semigroup such that almost all sample paths are cadlag.

\medskip

This gives a map $\Phi : \Omega \to D_X[0, \infty)$
sending $\omega \in \Omega$ to the sample path $t \mapsto B_t(\omega)$ (the map is defined $\mathbb{P}$-a.e.), and we thus obtain a
probability measure $\mathbb{P}_x$ on $(D_X[0, \infty), \mathcal{D})$ defined by $\mathbb{P}_x := (\Phi)_* \mathbb{P}$.
For each $t \geq 0$, we then obtain a probability measure $P(t,x, .)$ on $X$ defined by $P(t, x, .) := (\pi_t)_* \mathbb{P}_x$.

\medskip

Since $\pi_t \circ \Phi = B_t$, it is clear that the measure $P(t, x, .)$ on $X$
is just the distribution of $B_t$ for the process corresponding to the semigroup with initial distribution $\delta_x$.
The measures $P(t, x, .)$ are called the {\it transition probabilities}
associated to the semigroup $(T(t))_{t \geq 0}$.

\medskip

Since the process $(B_t)_{t \geq 0}$ corresponds to the semigroup $(T(t))_{t \geq 0}$,
and $B_0$ has initial distribution $\delta_x$, in fact we have for any $f \in B(X)$ and $t \geq 0$,
using (\ref{corresponds}) and the definition of conditional expectation,

\begin{align} \label{trans prob}
 \int_{\Omega} \mathbb{E}(f(B_t)|\mathcal{F}_0) d\mathbb{P} & = \int_{\Omega} (T(t)f)(B_0) d\mathbb{P} \nonumber \\
\Rightarrow \int_{\Omega} f(B_t) d\mathbb{P} & = (T(t)f)(x) \nonumber \\
\Rightarrow \int_X f(y) dP(t,x,y) & = (T(t)f)(x)
\end{align}
(where in the last line we used $(B_t)_* \mathbb{P} = P(t, x, .)$).
Thus the action of the semigroup $(T(t))_{t \geq 0}$ on functions is given by integrating against the transition probabilities.

\medskip

For $t \geq 0$, let $\theta_t : D_X[0, \infty) \to D_X[0, \infty)$ be the shift by time $t$ map, defined by $\theta_t(\gamma)(s) = \gamma(t+s)$.
Then it can be shown that for any $x \in X$ and $t \geq 0$, the measure $(\theta_t)_* \mathbb{P}_x$ is given by a convex combination of the
measures $\mathbb{P}_y, y \in Y$, namely
\begin{equation} \label{convex comb}
(\theta_t)_* \mathbb{P}_x = \int_X \mathbb{P}_y \ dP(t,x,y)
\end{equation}
(meaning that both sides above agree when applied to any $E \in \mathcal{D}$).

\medskip

Finally, if we assume that for any $x \in X$ and $\epsilon > 0$, the transition probabilities satisfy
\begin{equation} \label{small prob}
\lim_{t \to 0} \frac{1}{t} P(t, x, X - B(x, \epsilon)) = 0
\end{equation}
then for any probability measure $\nu$ on $X$, there is a Markov process $(B_t)_{t \geq 0}$ with initial distribution $\nu$, which
corresponds to the semigroup, such that almost all sample paths are continuous (\cite{ethierkurtz}, Proposition 2.9 of Chapter 4).
Thus with hypothesis (\ref{small prob}) above,
we obtain in this case for each $x \in X$ a probability measure $\mathbb{P}_x$
on $(C_X[0, \infty), \mathcal{C})$, where $C_X[0, \infty)$ is the space of continuous paths on $X$ and $\mathcal{C}$ is the sigma-algebra
generated by the coordinate maps $\pi_t : C_X[0, \infty) \to X, t \geq 0$ (it can be shown that this sigma-algebra coincides with the
Borel sigma-algebra of $C_X[0, \infty)$ when $C_X[0, \infty)$ is equipped with the topology of uniform convergence on compact sets).
The same relation (\ref{convex comb}) holds for the measures $\mathbb{P}_x$ on $C_X[0, \infty)$.

\medskip

\subsection{The heat semigroup and Brownian motion}

\medskip

Now let $X$ be a complete, simply connected Riemannian manifold of pinched negative sectional curvature. We indicate briefly how the
heat semigroup gives a Feller semigroup on $B(X)$, whose transition probabilities satisfy the hypothesis (\ref{small prob}). We will
obtain therefore for each $x \in X$, a Markov process $(B_t)_{t \geq 0}$ with initial distribution $\delta_x$, which
corresponds to the heat semigroup, such that almost all sample paths are continuous. This process is called the {\it Brownian motion}
on $X$ started at $x$. We obtain also a probability measure $\mathbb{P}_x$ on the space $C_X[0, \infty)$ of continuous paths in $X$,
called the {\it Wiener measure} on paths started at $x$.

\medskip

The heat kernel on any complete Riemannian manifold $X$ satisfies
$$
\int_X p(t,x,y) dvol(y) \leq 1
$$
for all $t > 0, x \in X$ (see section 2, Chapter VIII \cite{chavel}), so it follows that $(e^{t\Delta}f)(x) = \int_X p(t,x,y)f(y) dvol(y)$ is
well-defined for any $f \in B(X)$ and satisfies $||e^{t\Delta}f||_{\infty} \leq ||f||_{\infty}$. Thus the heat semigroup defines
a contraction semigroup on $B(X)$. Positivity of the semigroup is clear since the heat kernel is positive.

\medskip

We check the Feller property (c). Given $f \in C_0(X)$ and $\epsilon > 0$, choose a ball $B(x_0, R)$ such that
$|f(y)| \leq \epsilon/2$ for $y \in X - B(x_0, R)$ and choose $M > 2||f||_{\infty}$.
Given $t > 0$, it follows from Lemma \ref{outside ball} that
there is a constant $r_t > 0$ such that
$$
\int_{X - B(x, r_t)} p(t,x,y) dvol(y) \leq \epsilon/M
$$
for all $x \in X$. Then for $x$ such that $d(x,x_0) > R + r_t$, we have $B(x, r_t) \subset X - B(x_0, R)$, thus
\begin{align*}
\left|\int_X p(t,x,y)f(y) dvol(y)\right| & \leq \int_{X - B(x, r_t)} p(t,x,y)|f(y)| dvol(y) + \int_{B(x, r_t)} p(t,x,y)|f(y)| dvol(y) \\
                            & \leq \frac{M}{2} \frac{\epsilon}{M} + \frac{\epsilon}{2} \int_{B(x, r_t)} p(t,x,y) dvol(y) \\
                            & \leq \frac{\epsilon}{2} + \frac{\epsilon}{2} = \epsilon.
\end{align*}
It follows that $e^{t\Delta}f \in C_0(X)$.

\medskip

The condition that the heat semigroup be conservative, $e^{t \Delta} 1 = 1$, is equivalent to the condition
$$
\int_X p(t,x,y) dvol(y) = 1
$$
for all $x \in X, t > 0$. This is also referred to as {\it stochastic completeness} of $X$. Yau showed that any complete Riemannian manifold
with Ricci curvature bounded below is stochastically complete \cite{yau}. This holds in our case since the sectional curvature of $X$
is bounded below.

\medskip

For the strong continuity of the heat semigroup we use the fact that any $f \in C_0(X)$ is uniformly continuous. Given $\epsilon > 0$ let
$\delta > 0$ be such that $|f(x) - f(y)| < \epsilon/2$ for $d(x,y) < \delta$. From Lemma \ref{outside ball} it follows that there is a $t_0 > 0$
such that for $0 < t < t_0$ we have
$$
\int_{X - B(x, \delta)} p(t,x,y) dvol(y) \leq \frac{\epsilon}{M}
$$
for all $x \in X$ (where $M$ is a constant such that $M > 4||f||_{\infty}$). Then for any $x \in X$, for $0 < t < t_0$ we have
\begin{align*}
|(e^{t\Delta}f)(x) - f(x)| & = \left|\int_X p(t,x,y)(f(y) - f(x)) dvol(y)\right| \\
                           & \leq \int_{B(x, \delta)} p(t,x,y)|f(x) - f(y)| dvol(y) + \int_{X - B(x, \delta)} p(t,x,y)|f(x) - f(y)| dvol(y) \\
                           & \leq \frac{\epsilon}{2} \int_{B(x, \delta)} p(t,x,y) dvol(y) + 2||f||_{\infty} \frac{\epsilon}{M} \\
                           & \leq \frac{\epsilon}{2} + \frac{\epsilon}{2} = \epsilon
\end{align*}
Thus $||e^{t\Delta}f - f||_{\infty} \leq \epsilon$ for $0 < t < t_0$, so the heat semigroup is strongly continuous.

\medskip

This establishes that the heat semigroup is a Feller semigroup. It remains to check that the transition probabilities
$P(t,x,.)$ for the heat semigroup satisfy the condition (\ref{small prob}). From (\ref{trans prob}), the transition probability
$P(t,x,.)$ is a measure on $X$ satisfying
$$
\int_X f(y) dP(t,x,y) = (e^{t\Delta}f)(x)
$$
for all continuous bounded functions $f$ on $X$. It follows that the measure $P(t,x,.)$ is given by $dP(t,x,y) = p(t,x,y) dvol(y)$.
Now condition (\ref{small prob}) becomes
$$
\lim_{t \to 0} \frac{1}{t} \int_{X - B(x, \epsilon)} p(t,x,y) dvol(y) = 0
$$
which follows immediately from Lemma \ref{outside ball}. This finishes the construction of Brownian
motion on $X$.

\medskip

It is customary to write $\mathbb{E}_x(\phi)$ for the expectation of a function $\phi$ on the sample space $\Omega$
of a Brownian motion started at $x \in X$. We then have the following fundamental formula: for any $t > 0$ and $x \in X$,
$$
\mathbb{E}_x(f(B_t)) = (e^{t\Delta}f)(x)
$$
for any $f \in B(X)$. In particular, for any Borel set $E \subset X$, letting $f = 1_E$ gives
$$
\mathbb{P}_x( B_t \in E) = \int_E p(t,x,y) dvol(y)
$$
so that $p(t,x,.)$ is the density of the distribution of $B_t$ for Brownian motion started at $x$.

\medskip

Finally, it will be useful to note that the sample space $\Omega$ of Brownian motion on $X$ can always be taken to be
$C_X[0, \infty)$ with the random variables $B_t$ given by the coordinate maps $\pi_t : C_X[0, \infty) \to X$. Indeed, given a
Brownian motion $(B_t)_{t \geq 0}$ defined on some probability space $(\Omega, \mathcal{F}, \mathbb{Q})$, we can define the map
$\Phi : \Omega \to C_X[0, \infty)$ sending $\omega \in \Omega$ to the sample path $t \mapsto B_t(\omega)$, and let $\mathbb{P}$ be the
probability measure on $C_X[0, \infty)$ defined by $\mathbb{P} = (\Phi)_* \mathbb{Q}$. Then the $X$-valued process $(\pi_t)_{t \geq 0}$
on the probability space $(C_X[0, \infty), \mathcal{C}, \mathbb{P})$ is a Brownian motion whose sample paths have the same distribution
as those of $(B_t)_{t \geq 0}$. This realization of Brownian motion will be referred to as the {\it canonical coordinate process}.

\medskip

\section{The infinitesimal generator of the heat semigroup}

\medskip

For a Markov process $(B_t)_{t \geq 0}$ on $X$
corresponding to a semigroup $(T(t))_{t \geq 0}$ on $B(X)$, the {\it infinitesimal generator} of the semigroup is the
operator $A$ with domain $D(A) \subset B(X)$ defined by
$$
Af := \lim_{t \to 0} \frac{T(t)f - f}{t}
$$
where the domain $D(A)$ is the set of $f \in B(X)$ for which the above limit exists in $(B(X), ||.||_{\infty})$.
It is well known that if the semigroup is strongly continuous on $C_0(X)$ then the domain $D(A)$ is dense in $C_0(X)$
and the operator $A$ is closed.

\medskip

For the heat semigroup $(e^{t\Delta})_{t \geq 0}$, it is natural to expect that the infinitesimal generator
should be the Laplacian $\Delta$ acting on a suitable space of functions. We show the following:

\medskip

\begin{prop} \label{generator} For any $\phi \in C^{\infty}_c(X)$,
$$
\lim_{t \to 0} \frac{e^{t\Delta}\phi - \phi}{t} = \Delta \phi \ \ \hbox{ in } \ (B(X), ||.||_{\infty})
$$
Thus $C^{\infty}_c(X)$ is contained in the domain of the infinitesimal generator of the heat semigroup,
and the generator restricted to $C^{\infty}_c(X)$ equals the Laplacian $\Delta$.
\end{prop}

\medskip

\noindent{\bf Proof:} Recall Green's identity: for any $f, g \in C^{\infty}(X)$ and any precompact domain $D \subset X$
with smooth boundary, we have
$$
\int_D (f \Delta g - g \Delta) dvol = \int_{\partial D} \left(f \frac{\partial g}{\partial n} - g \frac{\partial f}{\partial n}\right) dA \ ,
$$
where $\frac{\partial f}{\partial n}, \frac{\partial g}{\partial n}$ are the normal derivatives on $\partial D$ and $dA$ is the Riemannian
volume measure of $\partial D$.

\medskip

For $f \in C^{\infty}(X)$ and $\phi \in C^{\infty}_c(X)$, applying the above formula on a domain $D$ such that
$supp(\phi) \subset D$ (so that $\phi = \frac{\partial \phi}{\partial n} = 0$ on $\partial D$) gives
\begin{equation} \label{symmetric}
\int_X f \Delta \phi \ dvol = \int_X f \Delta \phi \ dvol
\end{equation}

Given $\phi \in C^{\infty}_c(X)$ and $x \in X$, let $u(t) = (e^{t\Delta}\phi)(x), t \geq 0$. Then for $t > 0$,
$u(t) = \int_X p(t,x,y)\phi(y) dvol(y)$, so for $t > 0$ using (\ref{symmetric}) above we have
\begin{align*}
u'(t) & = \int_X \frac{\partial p}{\partial t}(t,x,y)\phi(y) dvol(y) \\
      & = \int_X \Delta_y p(t,x,y) \phi(y) \\
      & = \int_X p(t,x,y) \Delta \phi(y) dvol(y) \\
      & \to \Delta \phi(x) \ \ \hbox{ as } \ t \to 0 .
\end{align*}
Since $\lim_{t \to 0} u'(t)$ exists, it follows that $u$ is $C^1$ on $[0, \infty)$, and we can write
$$
e^{t\Delta}\phi(x) - \phi(x) = u(t) - u(0) = \int_{0}^{t} u'(s) ds = \int_{0}^{t} (e^{s\Delta} \Delta \phi)(x) ds
$$
It follows that for $t > 0$,
\begin{align*}
\left\|\frac{e^{t\Delta}\phi - \phi}{t} - \Delta \phi \right\|_{\infty} & \leq \frac{1}{t} \int_{0}^{t} || e^{s\Delta} \Delta \phi - \Delta \phi ||_{\infty} ds \\
                                                                        & \leq \sup_{0 \leq s \leq t} || e^{s\Delta} \Delta \phi - \Delta \phi ||_{\infty} \\
                                                                        & \to 0 \ \ \hbox{ as } \ t \to 0
\end{align*}
using the strong continuity of the heat semigroup. $\diamond$

\medskip

\begin{cor} \label{distributions} For $h \in B(X)$ with compact support,
$$
\frac{e^{t\Delta}h - h}{t} \to \Delta h
$$
in the sense of distributions on $X$.
\end{cor}

\medskip

\noindent{\bf Proof:} Let $\mathcal{D}(X)$ be the space of distributions on $X$ and let $(.,.)$ denote the pairing between
$\mathcal{D}(X)$ and $C^{\infty}_c(X)$. Then for $h \in B(X)$ with compact support and $\phi \in C^{\infty}_c(X)$, the distribution $\Delta h$ is
defined by $( \Delta h, \phi) := (h, \Delta \phi)$. We have,
using self-adjointness of the operators $e^{t\Delta}$ on $L^2(X)$, that
\begin{align*}
\left(\frac{e^{t\Delta} - id}{t}h, \phi \right) & = \left( h , \frac{e^{t\Delta} - id}{t} \phi \right) \\
                                                & \to (h, \Delta \phi) = (\Delta h, \phi) \ \ \hbox{ as } \ t \to 0
\end{align*}
since $\frac{e^{t\Delta} - id}{t} \phi \to \Delta \phi$ with respect to $||.||_{\infty}$ by Proposition \ref{generator},
and $h \in L^1(X)$. $\diamond$

\medskip

\section{The Dirichlet heat kernel of a precompact domain}

\medskip

Let $D \subset X$ be a precompact domain in $X$ with smooth boundary $\partial D$. Let $H^1(D)$ be the Sobolev space
$$
H^1(D) := \{ f \in L^2(D) | |\nabla f| \in L^2(D) \}
$$
(where $\nabla f$ is understood in the sense of distributions), equipped with the Sobolev norm
$$
||f||^2_{H^1} := \int_D |f|^2 dvol + \int_D |\nabla f|^2 dvol
$$
Let $H^1_0(D) \subset H^1(D)$ be the closure of $C^{\infty}_c(D)$ with respect to the Sobolev norm.
The Laplacian $\Delta$ with domain
$C^{\infty}_c(D) \subset L^2(D, dvol)$ turns out to be essentially self-adjoint on $L^2(D)$. Its closure is called the {\it Dirichlet Laplacian}
on $D$ and is denoted by $\Delta_D$. The domain of the Dirichlet Laplacian is given by
$$
Dom(\Delta_D) = \{ f \in H^1_0(D) | \Delta f \in L^2(D) \}
$$
(where $\Delta f$ is understood in the sense of distributions).

\medskip

In this case $L^2(D)$ has an orthonormal basis of
smooth eigenfunctions $\{ \phi_j \}_{j \geq 1}$ of the Dirichlet Laplacian (\cite{chavel}) with eigenvalues
$$
0 > \lambda_1(D) > \lambda_2(D) \geq \lambda_3(D) \geq ... \geq \lambda_k(D) \to -\infty \ \hbox{ as } k \to \infty
$$
(the multiplicity of the eigenvalue $\lambda_1(D)$ is one). The {\it Dirichlet heat semigroup} $(e^{t\Delta_D})_{t \geq 0}$ can
be defined using the functional calculus for self-adjoint operators, and it admits an integral kernel $p_D(t,x,y)$, called
the {\it Dirichlet heat kernel}. The Dirichlet heat kernel has the following expansion with respect to the basis $\{ \phi_j \}$:
$$
p_D(t,x,y) = \sum_{k = 1}^{\infty} e^{\lambda_k t} \phi_k(x) \phi_k(y)
$$
(the series converges absolutely and uniformly on $\overline{D}$ for all $t > 0$).

\medskip

The Dirichlet heat kernel can also be obtained as the smallest positive fundamental solution of the heat equation in $D$
with Dirichlet boundary conditions. The kernel $p_D : (0, \infty) \times \overline{D} \times \overline{D}$ is positive,
continuous, smooth on $(0, \infty) \times D \times D$, and satisfies
\begin{align*}
\frac{\partial p_D}{\partial t} & = \Delta_y p_D \\
p_D(t,x,.) & \to \delta_x \ \hbox{ as } t \to 0 \\
p_D(t,x,y) & = 0 \ \hbox{for } y \in \partial D
\end{align*}
(where the second condition above means $\int_D p_D(t,x,y)f(y) dvol(y) \to f(x)$ as $t \to 0$ for all continuous bounded functions $f$ on $D$).

\medskip

The Dirichlet heat kernel $p_D(t,x,y)$ can be obtained from the heat kernel $p(t,x,y)$ as follows (see \cite{chavel}, Chapter VII):

\medskip

As shown in \cite{chavel}, Chapter VII, for each $x \in \overline{D}$, there is a solution to the heat equation which is a
continuous function $g(., x, .) : [0, \infty) \times \overline{D} \to \mathbb{R}$
such that $g$ is $C^{\infty}$ on $(0, \infty) \times D$, satisfying
\begin{align*}
\frac{\partial g}{\partial t} & = \Delta_y g \ \ \hbox{ in } (0, \infty) \times D \\
g(0,x,y) & = 0 \ \ \hbox{ for all } y \in \overline{D} \\
g(t,x,y) & = -p(t,x,y) \ \ \hbox{ for all } t > 0, y \in \partial D
\end{align*}
The Dirichlet heat kernel is then given by
\begin{equation} \label{pd from p}
p_D(t,x,y) = p(t, x, y) + g(t, x, y) \ \ \hbox{ for } t > 0, x,y \in \overline{D}
\end{equation}

To proceed further, we will need the following well-known {\it parabolic maximum principle} for the heat equation (see \cite{chavel}, section VIII.1):

\medskip

For $T > 0$, let $D_T = (0, T) \times D$, and define the {\it parabolic boundary} $\partial_P D_T$ of $D_T$ by
$\partial_P D_T := (\{ 0 \} \times \overline{D}) \cup ([0, T] \times \partial D)$. Then we have:

\medskip

\begin{prop} \label{para max principle} (Parabolic maximum principle). Let $u$ be a continuous function on $\overline{D_T}$
which is a $C^{\infty}$ solution of the heat equation on $D_T$. Then
$$
\sup_{[0, T] \times \overline{D}} u = \sup_{\partial_P D_T} u
$$
and
$$
\inf_{[0, T] \times \overline{D}} u = \inf_{\partial_P D_T} u \ .
$$
\end{prop}

\medskip

It follows from the parabolic maximum principle that the function $g$ defined above satisfies $g(t,x,y) \leq 0$ for all $t \geq 0, x,y \in \overline{D}$, and hence
\begin{equation} \label{pd less than p}
p_D(t,x,y) \leq p(t,x,y) \ , \ t > 0, x,y \in \overline{D}.
\end{equation}
where $p(t,x,y)$ is as before the heat kernel of the whole manifold $X$. In particular,
$$
\int_D p_D(t,x,y) dvol(y) \leq 1,
$$
and the estimates (\ref{short time}), (\ref{long time}) also apply to $p_D(t,x,y)$ for $x,y \in \overline{D}$.

\medskip

Moreover, away from the boundary of $D$ and for small times we have the following estimate for the function $g$ which will be useful:

\medskip

\begin{lemma} \label{g small} There are constants $\alpha, \beta > 0$, such that for any compact $K \subset D$, if $\delta > 0$ denotes the
distance from $K$ to the boundary of $D$, then for all $t \in (0,1), x \in K$ we have
$$
\sup_{y \in \overline{D}} (-g(t,x,y)) \leq \alpha e^{-\delta^2 \beta/t}
$$
\end{lemma}

\medskip

\noindent{\bf Proof:} Given $x \in K$ and $t \in (0,1)$, since $g(0, x, y) = 0$ for all $y \in \overline{D}$, the parabolic
maximum principle gives
\begin{align*}
\sup_{y \in \overline{D}} (-g(t,x,y)) & = \sup_{(u,y) \in \partial_P D_t} (-g(u,x,y)) \\
                                      & \leq \sup_{ 0 < u \leq t, y \in \partial D} p(u,x,y)
\end{align*}
and the required estimate now follows from the estimate (\ref{short time}) for the heat kernel after choosing $\alpha > 0$ large enough and $\beta > 0$ small
enough and noting that $d(x,y) \geq \delta$ for $x \in K, y \in \partial D$. $\diamond$

\medskip

\section{Exit times and eigenfunctions}

\medskip

Let $(B_t)_{t \geq 0}$ be a Brownian motion on $X$ started at $x \in X$,
realized as the canonical coordinate process $(\pi_t)_{t \geq 0}$ on
$(C_X[0, \infty), \mathcal{C}, \mathbb{P}_x)$. For $t \geq 0$, we let $\theta_t : C_X[0, \infty) \to C_X[0, \infty)$
be the shift by time $t$ map as before. For convenience, we will denote $C_X[0, \infty)$ by $\Omega$ in what follows.

\medskip

We say that the Brownian motion is {\it transient} if
$$
\mathbb{P}_x( d(x, B_t) \to \infty \ \hbox{ as } \ t \to \infty) = 1 .
$$
It is well-known (see for eg. \cite{grigoryan2}, section 5) that the Brownian motion is transient if and only if
$$
\int_{1}^{\infty} p(t,x,x) dt < \infty .
$$
In our case we have the estimate (\ref{long time}), from which it is clear that the above integral converges, since $\lambda_1 < 0$.
Thus the Brownian motion is transient, which means with probability one it leaves every compact set in $X$ eventually.

\medskip

Given a precompact domain $D \subset X$, the {\it exit time} from $D$,
$\tau = \tau_D : C_X[0, \infty) \to [0, \infty]$, is defined by
$$
\tau(\gamma) := \inf \{ t > 0 : \gamma(t) \in X - D \} \ .
$$
Since the Brownian motion is transient, $\mathbb{P}_x( \tau < \infty) = 1$. So
the {\it exit point} from $D$, $\pi = \pi_D : \{ \tau < \infty \} \subset C_X[0, \infty) \to X$, given by
$$
\pi(\gamma) := \gamma(\tau(\gamma)) \ ,
$$
is defined almost everywhere. We note that if the starting point $x$ lies in $D$, then by continuity of the sample paths
the exit point lies on $\partial D$ almost surely, i.e. $\mathbb{P}_x(\pi \in \partial D) = 1$. The random variable 
$\pi$ is commonly written as $B_{\tau}$. 

\medskip

The following relation between the exit time $\tau$ and the Dirichlet heat kernel $p_D(t,x,y)$ is well-known: for
any $t > 0$ and $x \in D$, we have
\begin{equation} \label{exit distn}
\mathbb{P}_x(\tau \geq t) = \int_D p_D(t,x,y) dvol(y)
\end{equation}

\medskip

This leads to the following proposition:

\medskip

\begin{prop} \label{finite measure} For any $\lambda \in \mathbb{C}$ in the half-plane $\{ \hbox{Re} \ \lambda > \hbox{Re} \ \lambda_1 \}$,
there is a constant $C_{\lambda}$ such that
$$
\int_{\Omega} |e^{-\lambda \tau} | d\mathbb{P}_x \leq C_{\lambda}
$$
for all $x \in D$. Thus the complex measure $e^{-\lambda \tau} d\mathbb{P}_x$ has finite total variation for all $x \in D$.
\end{prop}

\medskip

\noindent{\bf Proof:} Since $|e^{-\lambda s}| = e^{-(\hbox{Re } \lambda)s}$ for $s$ real and $\mathbb{P}_x$ is a probability
measure, if $\hbox{Re } \lambda \geq 0$ we are done, so
we may as well assume that $\lambda$ is real and $\lambda_1 < \lambda < 0$. Let $G$ be the
 monotone decreasing function $G(s) = \mathbb{P}_x( \tau \geq s)$. From the relation (\ref{exit distn}), using $p_D(t,x,y) \leq p(t,x,y)$
 and the estimate (\ref{long time}) it follows that $G(s) \leq C_1 e^{\lambda_1 s}$ for some constant $C_1 > 0$
 independent of $x$, so $e^{-\lambda s} G(s) \to 0$ as $s \to \infty$.
 We can then integrate by parts to obtain:
\begin{align*}
\int_{\Omega} e^{-\lambda \tau} d\mathbb{P}_x & = - \int_{0}^{\infty} e^{-\lambda s} dG(s) \\
                                              & = - [ e^{-\lambda s} G(s) ]_{0}^{\infty} + \int_{0}^{\infty} (-\lambda) e^{-\lambda s} G(s) ds \\
                                              & \leq 1 + \int_{0}^{1} (-\lambda)e^{-\lambda s} ds + C_1 \int_{1}^{\infty} e^{(\lambda_1 - \lambda)s} ds \\
                                              & = C_{\lambda}
\end{align*}
$\diamond$

\medskip

We observe that for any $t > 0$ we have
\begin{equation} \label{shift exit}
\tau \circ \theta_t = \tau - t \ \ \hbox{ on the set } \ \{  t \leq \tau < \infty \} \ .
\end{equation}

\medskip

We will need the following estimate on the $\mathbb{P}_x$-measure of the set $\{ \tau < t \}$ for $x$ in a compact in $D$ and for small times $t$:

\medskip

\begin{lemma} \label{exit time large small prob} There are constants $\gamma, \beta' > 0$, such that for any compact $K \subset D$,
if $\delta > 0$ denotes the distance from $K$ to the boundary of $D$, then for all $x \in K$ and $0 < t < min(1, \delta/(Dh))$ we have
$$
\mathbb{P}_x(\tau < t) \leq \gamma e^{-\delta^2 \beta' /t} \ .
$$
\end{lemma}

\noindent{\bf Proof:} Given $x \in K$ and $t \in (0,1)$, from (\ref{exit distn}) and using (\ref{pd from p}),
Lemma \ref{g small} and Lemma \ref{outside ball} we have
\begin{align*}
\mathbb{P}_x(\tau < t) & = 1 - \int_{D} p_D(t,x,y) dvol(y) \\
                       & = (1 - \int_{D} p(t,x,y) dvol(y)) + \int_{D} (-g(t,x,y)) dvol(y) \\
                       & \leq \int_{X - B(x, \delta)} p(t,x,y) dvol(y) + vol(D) \alpha e^{-\delta^2 \beta/t} \\
                       & \leq \kappa e^{-\delta^2 \eta/t} + vol(D) \alpha e^{-\delta^2 \beta/t} \\
                       & \leq  \gamma e^{-\delta^2 \beta' /t}
\end{align*}
where $\beta' = min(\eta, \beta)$ and $\gamma = \kappa + \alpha vol(D)$. $\diamond$

\medskip

We now come to the proof of Theorem \ref{mainthm}. We fix a $\lambda \in \mathbb{C}$ with $\hbox{Re}(\lambda) > \lambda_1$.
For $x \in D$, we define the measure $\mu_{x, \lambda}$ on $\partial D$ by
$$
\mu_{x, \lambda} := \pi_* ( e^{-\lambda \tau} d\mathbb{P}_x ) \ ,
$$
where $\pi$ is the exit point map defined previously. Since $\mathbb{P}_x (\pi \in \partial D) = 1$ for $x \in D$, the measure $\mu_{x, \lambda}$
is supported on $\partial D$. By Proposition \ref{finite measure}, $\mu_{x, \lambda}$ is a complex measure on $\partial D$ of finite total variation.

\medskip

For the rest of this section, we fix a continuous function $\phi : \partial D \to \mathbb{C}$.
We then define a function $h$ on $D$ by
$$
h(x) := \int_{\partial D} \phi(y) d\mu_{x, \lambda}(y) \ .
$$
We wish to show that $h$ is an eigenfunction of $\Delta$ on $D$ with eigenvalue $\lambda$, and that $h(x) \to \phi(\xi)$ as $x \to \xi \in \partial D$.

\medskip

It will be convenient to write the function $h$ as follows: we first extend $\phi$ to a function on all of $X$ such that $\phi$ is continuous
with compact support. We then define a function $\Phi : \{ \tau < \infty \} \subset C_X[0, \infty) \to \mathbb{C}$ by
$$
\Phi(\gamma) := \phi(\pi(\gamma)) \ .
$$
We will also write $\phi(B_{\tau})$ for the random variable $\Phi$. From the definition of $\mu_{x, \lambda}$, we can write $h$ as
\begin{equation} \label{h defn}
h(x) = \int_{\Omega} e^{-\lambda \tau} \Phi d\mathbb{P}_x = \int_{\Omega} e^{-\lambda \tau} \phi(B_{\tau}) d\mathbb{P}_x 
\end{equation}
We note that if $x \in X - \overline{D}$, then $\mathbb{P}_x( \tau = 0) = 1$ and $\mathbb{P}_x( \pi = x) = 1$,
so if we define $h$ by equation (\ref{h defn}) above, then $h(x) = \phi(x)$ for $x \in X - \overline{D}$. Thus we can regard $h$ as
a bounded function on $X$ with compact support, defined by (\ref{h defn}) for all $x \in X$.

\medskip

\begin{lemma} \label{h and semigroup} For any compact $K \subset D$, we have
$$
\frac{e^{t\Delta}h - h}{t} \to \lambda h
$$
uniformly on $K$ as $t \to 0$.
\end{lemma}

\medskip

\noindent{\bf Proof:} We first note that for $t > 0$, $\pi \circ \theta_t = \pi$ on the set $\{ t \leq \tau < \infty \}$, and
hence $\Phi \circ \theta_t = \Phi$ on $\{ t \leq \tau < \infty \}$. Together with the relations (\ref{shift exit}) and (\ref{convex comb}),
this leads to
the following, for $x \in K \subset D$:
\begin{align*}
(e^{t\Delta}h)(x) & = \int_X p(t,x,y) h(y) dvol(y) \\
                  & = \int_X p(t,x,y) (\int_{\Omega} e^{-\lambda \tau} \Phi d\mathbb{P}_y) dvol(y) \\
                  & = \int_{\Omega} e^{-\lambda \tau} \Phi d((\theta_t)_* \mathbb{P}_x) \\
                  & = \int_{\Omega} e^{-\lambda \tau \circ \theta_t} \Phi \circ \theta_t d\mathbb{P}_x \\
                  & = \int_{\{ t \leq \tau < \infty \}} e^{-\lambda \tau \circ \theta_t} \Phi \circ \theta_t d\mathbb{P}_x + \int_{\{\tau < t\}} e^{-\lambda \tau \circ \theta_t} \Phi \circ \theta_t d\mathbb{P}_x \\
                  & = \int_{\{t \leq \tau < \infty\}} e^{-\lambda (\tau - t)} \Phi d\mathbb{P}_x + \int_{\{\tau < t\}} e^{-\lambda \tau \circ \theta_t} \Phi \circ \theta_t d\mathbb{P}_x \\
                  & = \int_{\Omega} e^{-\lambda (\tau - t)} \Phi d\mathbb{P}_x - \int_{\{\tau < t\}} e^{-\lambda (\tau - t)} \Phi d\mathbb{P}_x + \int_{\{\tau < t\}} e^{-\lambda \tau \circ \theta_t} \Phi \circ \theta_t d\mathbb{P}_x \\
                  & = e^{\lambda t} h(x) - A(x,t) + B(x,t)
\end{align*}
where $A(x,t) = \int_{\{\tau < t\}} e^{-\lambda (\tau - t)} \Phi d\mathbb{P}_x$ and
$B(x,t) = \int_{\{\tau < t\}} e^{-\lambda \tau \circ \theta_t} \Phi \circ \theta_t d\mathbb{P}_x$. Thus
\begin{equation} \label{good plus errors}
\frac{(e^{t\Delta}h)(x) - h(x)}{t} = \frac{e^{\lambda t} - 1}{t} h(x) - \frac{1}{t}A(x,t) + \frac{1}{t}B(x,t) \ ,
\end{equation}
so to finish the proof it suffices to show that $A(x,t) = o(t), B(x,t) = o(t)$ as $t \to 0$, uniformly in $x \in K$.

\medskip

Let $\delta > 0$ be the distance from $K$ to the boundary of $D$, and let $M > 0$ be such that $|\phi| \leq M$ on $X$ (and so $|\Phi| \leq M$ on $\{ \tau < \infty \}$).
Then for $0 < t < min(1, \delta/(Dh))$ and $x \in K$, using Lemma \ref{exit time large small prob}
we have
\begin{align*}
\frac{1}{t}|A(x,t)| & \leq \frac{1}{t} \int_{\{\tau < t\}} |e^{-\lambda \tau} e^{\lambda t} \Phi| d\mathbb{P}_x \\
                    & \leq \frac{1}{t} e^{2|\lambda| t} M \mathbb{P}_x( \tau < t) \\
                    & \leq \frac{1}{t} e^{2|\lambda| t} M \gamma e^{-\delta^2 \beta' /t} \\
                    & \to 0 \ \ \hbox{ as } \ t \to 0, \ \hbox{ uniformly in } \ x \in K
\end{align*}
We now estimate the term $B(x,t)$. For this it will be convenient to use Holder's inequality.
Since $\hbox{Re}(\lambda) > \lambda_1$, we can chose $p > 1$ such that $\lambda' := p \lambda$ also satisfies $\hbox{Re}(\lambda') > \lambda_1$.
We can then estimate the $L^p$ norm of $e^{-\lambda \tau \circ \theta_t}$ with respect to $\mathbb{P}_x$ using the relation (\ref{convex comb})
and Proposition \ref{finite measure} as follows:
\begin{align*}
\int_{\Omega} |e^{-\lambda \tau \circ \theta_t}|^p d\mathbb{P}_x & = \int_X \left( \int_{\Omega} |e^{-\lambda' \tau}| d\mathbb{P}_y \right) p(t,x,y) dvol(y) \\
                                                                 & = \int_D \left( \int_{\Omega} |e^{-\lambda' \tau}| d\mathbb{P}_y \right) p(t,x,y) dvol(y)
                                                                 + \int_{X - \overline{D}} \left( \int_{\Omega} |e^{-\lambda' \tau}| d\mathbb{P}_y \right) p(t,x,y) dvol(y) \\
                                                                 & = \int_D \left( \int_{\Omega} |e^{-\lambda' \tau}| d\mathbb{P}_y \right) p(t,x,y) dvol(y)
                                                                 + \int_{X - \overline{D}} p(t,x,y) dvol(y) \\
                                                                 & \leq \int_D C_{\lambda'} p(t,x,y) dvol(y) + 1 \\
                                                                 & \leq C_{\lambda'} + 1
\end{align*}
(where we used $\mathbb{P}_y(\tau = 0) = 1$ for $y \in X - \overline{D}$ above). Letting $q \in (1, \infty)$ be such that $1/p + 1/q = 1$, we can estimate
$|B(x,t)|$ using Holder's inequality and Lemma \ref{exit time large small prob}:
\begin{align*}
\frac{1}{t}|B(x,t)| & \leq \frac{1}{t} \int_{\{\tau < t\}} |e^{-\lambda \tau \circ \theta_t}| |\Phi \circ \theta_t| d\mathbb{P}_x \\
                    & \leq \frac{1}{t} M \int_{\Omega} |e^{-\lambda \tau \circ \theta_t}| \cdot 1_{\{ \tau < t \}} d\mathbb{P}_x \\
                    & \leq \frac{M}{t} \left( \int_{\Omega} |e^{-\lambda \tau \circ \theta_t}|^p d\mathbb{P}_x \right)^{1/p} \mathbb{P}_x( \tau < t )^{1/q} \\
                    & \leq \frac{M}{t} ( C_{\lambda'} + 1 )^{1/p} \gamma e^{-\delta^2 \beta' / (qt)} \\
                    & \to 0 \ \ \hbox{ as } \ t \to 0, \ \hbox{ uniformly in } \ x \in K
\end{align*}
It now follows from (\ref{good plus errors}) that
$$
\frac{e^{t\Delta}h - h}{h} \to \lambda h
$$
uniformly on $K$ as $t \to 0$. $\diamond$

\medskip

We can now show that $h$ is an eigenfunction of $\Delta$ with eigenvalue $\lambda$:

\medskip

\begin{prop} \label{h eigenfunction} The function $h$ is $C^{\infty}$ on $D$ and satisfies
$$
\Delta h = \lambda h
$$
on $D$.
\end{prop}

\medskip

\noindent{\bf Proof:} Let $(.,.)$ denote the pairing between distributions on $D$ and functions in $C^{\infty}_c(D)$. Given by
$\psi \in C^{\infty}_c(D)$, let $K = supp(\psi) \subset D$. Then it follows from Lemma \ref{h and semigroup} above that
$$
\left( \frac{e^{t\Delta}h - h}{t}, \psi \right) \to ( \lambda h, \psi) \ \hbox{ as } \ t \to 0 \ .
$$
On the other hand, by Corollary \ref{distributions}, we have
$$
\left( \frac{e^{t\Delta}h - h}{t}, \psi \right) \to ( \Delta h, \psi) \ \hbox{ as } \ t \to 0 \ .
$$
It follows that $\Delta h = \lambda h$ as distributions on $D$, and hence by elliptic regularity $h$ is $C^{\infty}$ on $D$
and $\Delta h = \lambda h$ as functions on $D$. $\diamond$

\medskip

To complete the proof of Theorem \ref{mainthm}, we will need the following lemmas:

\medskip

\begin{lemma} \label{tau bigger t small prob}
Let $\xi \in \partial D$. Then for any $t > 0$,
\begin{equation} \label{p tau tends zero}
\mathbb{P}_x( \tau \geq t ) \to 0 \ \ \hbox{ as } \ x \in D \to \xi \in \partial D \ ,
\end{equation}
and
\begin{equation} \label{e tau tends zero}
\int_{ \{ \tau \geq t \} } |e^{-\lambda \tau}| d\mathbb{P}_x \to 0 \ \ \hbox{ as } \ x \in D \to \xi \in \partial D \ .
\end{equation}
\end{lemma}

\medskip

\noindent{\bf Proof:} Since for $t > 0$ fixed, $p_D(t,.,.)$ is a continuous function on $\overline{D} \times \overline{D}$
which vanishes for $(x,y) \in \partial D \times \overline{D}$, we have $p_D(t,x,y) \to 0$ uniformly in $y \in \overline{D}$
as $x \to \xi$, hence
$$
\mathbb{P}_x( \tau \geq t ) = \int_D p_D(t,x,y) dvol(y) \to 0 \ \ \hbox{ as } \ x \in D \to \xi \in \partial D \ ,
$$
which proves (\ref{p tau tends zero}) above. For (\ref{e tau tends zero}), since $|e^{-\lambda \tau}| \leq 1$ for $\hbox{Re } \ \lambda \geq 0$,
we may as well assume that $\lambda$ is real and $\lambda_1 < \lambda < 0$. For $x \in D$ we define the function
$$
G_x(s) := \mathbb{P}_x( \tau \geq s ) = \int_D p_D(s,x,y) dvol(y) \ , s > 0 \ ,
$$
then as before we have $G_x(s) \leq C_1 e^{\lambda_1 s}$ for some constant $C_1 > 0$ independent of $x$, so
\begin{align*}
\int_{ \{ \tau \geq t \} } e^{-\lambda \tau} d\mathbb{P}_x & = - \int_{t}^{\infty} e^{-\lambda s} dG_x(s) \\
                                                           & = e^{-\lambda t} \mathbb{P}_x( \tau \geq t ) + (-\lambda) \int_{t}^{\infty} e^{-\lambda s} G_x(s) ds \\
                                                           & = e^{-\lambda t} \mathbb{P}_x( \tau \geq t ) + (-\lambda) \int_D \int_{t}^{\infty} e^{-\lambda s} p_D(s,x,y) ds dvol(y)
\end{align*}
Now the first term on the right-hand side above tends to zero as $x \to \xi$ by (\ref{p tau tends zero}), while for the
second term we can use the dominated convergence theorem as follows: from estimate (\ref{long time}), we can find a constant $C_2 > 0$
such that for $x,y \in \overline{D}$ and $s \geq t$, we have $p(s,x,y) \leq C_2 e^{\lambda_1 s}$. This gives
\begin{align*}
\int_D \int_{t}^{\infty} e^{-\lambda s} p_D(s,x,y) ds dvol(y) & \leq \int_D \int_{t}^{\infty} e^{-\lambda s} p(s,x,y) ds dvol(y) \\
                                                              & \leq C_2 \int_D \int_{t}^{\infty} e^{-\lambda s} e^{\lambda_1 s} ds dvol(y) \\
                                                              & < +\infty
\end{align*}
since $\lambda > \lambda_1$. Hence dominated convergence applies, and so
$$
\int_D \int_{t}^{\infty} e^{-\lambda s} p_D(s,x,y) ds dvol(y) \to 0 \ \ \hbox{ as } \ x \to \xi
$$
since $p_D(s,x,y) \to 0$ as $x \to \xi$ for all $s \geq t, y \in D$. This proves (\ref{e tau tends zero}). $\diamond$

\medskip

For $x \in X$ and $\delta > 0$, we define the exit time $\tau_{x, \delta} : \Omega \to [0, \infty]$
of Brownian motion from the ball $B(x, \delta)$ by
$$
\tau_{x, \delta}(\gamma) := \inf \{ t > 0 : \gamma(t) \in X - B(x, \delta) \} \ .
$$

\begin{lemma} \label{p exit from ball}
For any $\delta > 0$,
$$
\mathbb{P}_x( \tau_{x, \delta} < t ) \to 0 \ \ \hbox{ as } \ t \to 0
$$
uniformly in $x \in X$.
\end{lemma}

\medskip

\noindent{\bf Proof:} Given $x \in X$, let $B = B(x, \delta)$ and let $p_B(.,.,.)$ be the Dirichlet heat kernel
of the ball $B$. As in section 5, we can write $p_B(.,x,.) = p(.,x,.) + g(.,x,.)$ where $g : [0, \infty) \times \overline{B}$ is continuous,
is a solution of the heat equation on $(0, \infty) \times B$, and satisfies the boundary conditions $g(0,x,.) \equiv 0$, $g(t,x,y) = - p(t,x,y)$ for
$t > 0, y \in \partial B$. Lemma \ref{g small} applies in this situation to give constants $\alpha, \beta > 0$ independent of $x$ such that
$$
\sup_{y \in \overline{B}} (-g(t,x,y)) \leq \alpha e^{-\delta^2 \beta/t}
$$
Together with Lemma \ref{outside ball}, this gives, for $0 < t < \delta/(Dh)$,
\begin{align*}
\mathbb{P}_x( \tau_{x, \delta} < t ) & = 1 - \int_B p_B(t, x, y) dvol(y) \\
                                     & = (1 - \int_B p(t,x,y) dvol(y)) + \int_B (-g(t,x,y)) dvol(y) \\
                                     & = \int_{X - B(x,\delta)} p(t,x,y) dvol(y) + \int_B (-g(t,x,y)) dvol(y) \\
                                     & \leq \kappa e^{-\delta^2 \eta/t} + vol(B(x,\delta)) \alpha e^{-\delta^2 \beta /t} \\
                                     & \to 0 \ \ \hbox{ as } \ t \to 0 \ \hbox{ uniformly in } \ x \in X
\end{align*}
(where we have used the fact that $vol(B(x, \delta)) \leq C_{\delta}$ for some constant independent of $x$, which holds
since the sectional curvature of $X$ is bounded below). $\diamond$

\medskip

We can now prove Theorem \ref{mainthm}:

\medskip

\noindent{\bf Proof of Theorem \ref{mainthm}:} It remains to show that $h(x) \to \phi(\xi)$ as $x \in D \to \xi \in \partial D$.
Let $\xi \in \partial D$, and fix $\epsilon > 0$. We choose $\delta > 0$ such that $|\phi(y) - \phi(\xi)| \leq \epsilon$
for $y \in \partial D$ with $d(y, \xi) \leq \delta$, and fix a constant $M > ||\phi||_{\infty}$. For any $t > 0$, we can write the space $\Omega$ as
$$
\Omega = \{ \tau < t, \tau < \tau_{\xi, \delta} \} \sqcup \{ \tau < t, \tau \geq \tau_{\xi, \delta} \} \sqcup \{ \tau \geq t \} \ ,
$$
where $\tau_{\xi, \delta}$ is the exit time from the ball $B(\xi, \delta)$ as defined previously.
We then have, for $x \in D$,
\begin{align} \label{c d e}
|h(x) - \phi(\xi)| & \leq \int_{\Omega} |e^{-\lambda \tau} \phi(B_{\tau}) - \phi(\xi)| d\mathbb{P}_x \nonumber \\
                   & = C(x,t) + D(x,t) + E(x,t) \ ,
\end{align}
where
\begin{align*}
C(x,t) & = \int_{ \{ \tau < t, \tau < \tau_{\xi, \delta} \} } |e^{-\lambda \tau} \phi(B_{\tau}) - \phi(\xi)| d\mathbb{P}_x \ , \\
D(x,t) & = \int_{ \{ \tau < t, \tau \geq \tau_{\xi, \delta} \} } |e^{-\lambda \tau} \phi(B_{\tau}) - \phi(\xi)| d\mathbb{P}_x \ , \\
E(x,t) & = \int_{ \{ \tau \geq t \} } |e^{-\lambda \tau} \phi(B_{\tau}) - \phi(\xi)| d\mathbb{P}_x \ .
\end{align*}
We will show that by first choosing $t > 0$ small enough, all the terms $C(x,t), D(x,t), E(x,t)$ are small for $x$ close enough to $\xi$.

\medskip

Note that on the set $\{ \tau < t, \tau < \tau_{\xi, \delta} \}$, for
$x \in B(\xi, \delta) \cap D$,we have, $\mathbb{P}_x$-a.s., that $B_{\tau} \in B(\xi, \delta) \cap \partial D$, and hence
$|\phi(B_{\tau}) - \phi(\xi)| \leq \epsilon$ on this set $\mathbb{P}_x$-a.s. for $x \in B(\xi, \delta) \cap D$.
Also on this set, since $\tau < t$, for $t > 0$ small enough
we have $|e^{-\lambda \tau} - 1| \leq 2|\lambda| t$. We can then estimate, for $x \in B(\xi, \delta) \cap D$,
\begin{align*}
C(x,t) & \leq \int_{ \{ \tau < t, \tau < \tau_{\xi, \delta} \} } |e^{-\lambda \tau}| |\phi(B_{\tau}) - \phi(\xi)| + |e^{-\lambda \tau} - 1| |\phi(\xi)| d\mathbb{P}_x \\
         & \leq e^{-\lambda_1 t} \cdot \epsilon + 2|\lambda|t \cdot M \\
         & \leq C\epsilon
\end{align*}
for all $t \in (0, t_1]$, for some constants $C, t_1 > 0$ independent of $x$.

\medskip

To estimate $D(x,t)$, we note that for $x \in B(\xi, \delta/2)$ we have $B(x, \delta/2) \subset B(\xi, \delta)$, and so
$\mathbb{P}_x( \tau_{x, \delta/2} \leq \tau_{\xi, \delta} ) = 1$ (since any path starting at $x$ must exit $B(x, \delta/2)$
before it exits $B(\xi, \delta)$), from which we get
\begin{align} \label{tau delta tau}
\mathbb{P}_x( \tau \geq \tau_{\xi, \delta} ) & = \mathbb{P}_x( \tau \geq \tau_{\xi, \delta}, \tau < t ) + \mathbb{P}_x( \tau \geq \tau_{\xi, \delta}, \tau \geq t ) \nonumber \\
                                             & \leq \mathbb{P}_x( \tau_{x, \delta/2} < t ) + \mathbb{P}_x( \tau \geq t )
\end{align}
for $x \in B(\xi, \delta/2)$. Now from Lemma \ref{p exit from ball}, it follows that there is $t_2 > 0$ such that for all $t \in (0, t_2]$ we have
$$
\mathbb{P}_x( \tau_{x, \delta/2} < t ) \leq \epsilon
$$
for all $x \in B(\xi, \delta/2)$. We may also assume choosing $t_2 > 0$ small enough
 that $t_2 \leq t_1$, and $e^{|\lambda| t_2} \leq 2$.
We now fix $t_2$. Then from Lemma \ref{tau bigger t small prob}, we can choose $\delta_1 \in (0, \delta/2)$ such that
$$
\mathbb{P}_x( \tau \geq t_2 ) \leq \epsilon
$$
for all $x \in B(\xi, \delta_1) \cap D$. With these choices, we get for all $x \in B(\xi, \delta_1) \cap D$ the estimate
\begin{align*}
D(x, t_2) & \leq \int_{ \{ \tau < t, \tau \geq \tau_{\xi, \delta} \} } |e^{-\lambda \tau} \phi(B_{\tau}) - \phi(\xi)| d\mathbb{P}_x \\
            & \leq (2 \cdot M + M) \mathbb{P}_x( \tau \geq \tau_{\xi, \delta} ) \\
            & \leq 3M \cdot (\epsilon + \epsilon)
\end{align*}
where we have used (\ref{tau delta tau}) and the way $t_2, \delta_1$ were chosen.

\medskip

Finally, we estimate $E(x, t_2)$:
\begin{align*}
E(x, t_2) & \leq \int_{ \{ \tau \geq t_2 \} } |e^{-\lambda \tau}\phi(B_{\tau})| + |\phi(\xi)| d\mathbb{P}_x \\
            & \leq M \int_{ \{ \tau \geq t_2 \} } |e^{-\lambda \tau}| d\mathbb{P}_x + M \mathbb{P}_x( \tau \geq t_2 ) \ .
\end{align*}
From Lemma \ref{tau bigger t small prob} it follows that there is a $\delta_2 \in (0, \delta_1)$
such that both terms on the right-hand side above are less than $\epsilon$ for $x \in B(\xi, \delta_2) \cap D$, so that $E(x,t_2) \leq 2 \epsilon$
for $x \in B(\xi, \delta_2) \cap D$.

\medskip

Putting together all the estimates, we get, for $x \in B(\xi, \delta_2) \cap D$,
\begin{align*}
|h(x) - \phi(\xi)| & \leq C(x, t_2) + D(x, t_2) + E(x, t_2) \\
                   & \leq (C + 6M + 2)\epsilon
\end{align*}
and so $h(x) \to \phi(\xi)$ as $x \to \xi$. $\diamond$

\bibliography{moeb}
\bibliographystyle{alpha}

\end{document}